\documentclass[graybox]{svmult}

\usepackage{type1cm}
\usepackage{makeidx}
\usepackage{graphicx}
\usepackage{multicol} 
\usepackage[bottom]{footmisc}
\usepackage{hyperref}
\usepackage{subcaption}
\usepackage{newtxtext} 
\usepackage[varvw]{newtxmath}
\makeindex

\begin{document}

\title*{Reduced order modeling for spectral element methods: current developments in Nektar++ and further perspectives}
\titlerunning{Reduced order modeling for SEM} 

\author{Martin W. Hess, Andrea Lario, Gianmarco Mengaldo and Gianluigi Rozza}

\institute{Martin W. Hess (Corresponding Author) \at SISSA mathLab, International School for Advanced Studies, via Bonomea 265, I-34136 Trieste, Italy, \email{ mhess@sissa.it}
\and Andrea Lario  \at SISSA mathLab, International School for Advanced Studies, via Bonomea 265, I-34136 Trieste, Italy, \email{ alario@sissa.it}
\and Gianmarco Mengaldo  \at NUS, National University of Singapore, 21 Lower Kent Ridge Rd, Singapore 119077, Singapore \email{mpegim@nus.edu.sg}
\and Gianluigi Rozza  \at SISSA mathLab, International School for Advanced Studies, via Bonomea 265, I-34136 Trieste, Italy, \email{grozza@sissa.it}
}

\maketitle

\abstract{In this paper, we present recent efforts to develop reduced order modeling (ROM) capabilities for spectral element methods (SEM). Namely, we detail the implementation of ROM for both continuous Galerkin and discontinuous Galerkin methods in the spectral/hp element library Nektar++. The ROM approaches adopted are intrusive methods based on the proper orthogonal decomposition (POD). They admit an offline-online decomposition, such that fast evaluations for parameter studies and many-queries are possible. An affine parameter dependency is exploited such that the reduced order model can be 
evaluated independent of the large-scale discretization size. 
The implementation in the context of SEM can be found in the open-source model reduction software ITHACA-SEM.
}

\section{Introduction}
\label{sec:1}
The spectral element method (SEM)~\cite{karniadakis1999spectral} is a discretization technique for partial differential equations (PDEs) that uses high-order polynomials on a tessellation of elements constituting the computational domain. 
Unlike other numerical discretizations, such as traditional finite element, finite difference and finite volume methods, SEM provides exponential convergence properties to the PDE solution and allows for tunable numerical properties, such as tunable diffusion and dispersion errors. 
This flexibility has been extensively leveraged to construct high-fidelity simulation capabilities~\cite{mengaldo2021industry} in the context of both continuous Galerkin (CG)~\cite{moura2016eigensolution,moura2020spatial}, and discontinuous Galerkin (DG) formulations~\cite{mengaldo2015discontinuous,moura2017eddy,moura2017setting,mengaldo2018spatial,moxey2017towards,mengaldo2018spatial,fernandez2019non,moura2020viscous,tonicello2021fully}. The former (CG) constructs a numerical discretization enforcing the solution to be continuous between elements. The latter (DG) enforces the numerical fluxes (and not the solution) to be continuous between elements. The use of SEM for high-fidelity simulations provides a pathway to next generation computational tools for engineering analysis and design, in fluid dynamics and other sectors. This offers the opportunity to devise improved reduced order models that may be used for real-time design and control purposes. 

In this paper, with present recent efforts that were undertaken to develop reduced order modeling (ROM)~\cite{SchildersMOR, hesthaven2015certified, deville_fischer_mund_2002, HessRozza2019, HessQuainiRozza2020} approaches within SEM, both for CG and DG methods. These efforts were conveyed in the implementation of ROM into the spectral/hp element library Nektar++\footnote{\texttt{www.nektar.info}}~\cite{cantwell2015nektar++,moxey2020nektar++}. In particular, we show the global assembly of the PDE solution and how the model reduction relates to the full-order solver for both CG and DG. The different solver modules of Nektar++ share basic functionalities such as the geometry description and quadrature rules, among others. This also holds true for the Nektar++ incompressible and compressible solver modules. The former uses a CG discretization, while the latter uses DG. As a consequence, the model reduction codes for incompressible and compressible flow simulations are independent from each other. The implementation of ROM methods in Nektar++ is available in the open-source model reduction software ITHACA-SEM\footnote{\texttt{https://mathlab.sissa.it/ITHACA-SEM}}.

The open-source software ITHACA-SEM currently has the capability to generate POD-based ROMs for 2D incompressible Navier-Stokes equations with parametric variation in geometry and/or kinematic viscosity. 
The parametric dependency on geometry parameters is assumed to be affine and the user can specify the affine form in a header file. After another compilation of ITHACA-SEM, the affine form is then available. It can thus serve as a head start for a developer seeking to work in ROMs with the SEM and also to a practitioner within the boundaries mentioned.
ITHACA-SEM is delivered and compiled with the Nektar++ master branch, which is periodically merged into the code.
A few test cases are part of the Nektar++ unit test and additional examples are available.

The paper is organized as follows. In section~\ref{sec:cg}, we present the CG method and the ROM approach implemented. In section~\ref{sec:dg}, we detail the DG method and the associated ROM method implemented. In section~\ref{sec:conclusions}, we draw some conclusions and future perspectives.




\section{Continuous Galerkin: incompressible flow simulations}\label{sec:cg}

\subsection{Overview}
Let $\Omega$ denote the spatial computational domain.
Incompressible, viscous fluid motion in the domain $\Omega$ over a time interval $(0, T)$ is governed by the incompressible  \emph{Navier-Stokes} equations
\begin{eqnarray}
\frac{\partial u}{\partial t} + u \cdot \nabla u &=& - \nabla p + \nu \Delta u + f, \label{Hess:NSE0} \\
\nabla \cdot u &=& 0.
\label{Hess:NSE1}
\end{eqnarray}
where $u$ is the fluid velocity, $p$ is the pressure, $\nu$ is the kinematic viscosity, and $f$ is a body forcing term. The following boundary and initial conditions 
\begin{eqnarray}
u &=& d \quad \text{ on } \Gamma_D \times (0, T), \\
\nabla u \cdot n &=& g \quad \text{ on } \Gamma_N \times (0, T), \\
u &=& u_0 \quad \text{ in } \Omega \times 0,
\label{Hess:NSE_boundaryCond}
\end{eqnarray}
fully define the problem~\eqref{Hess:NSE0}-\eqref{Hess:NSE1}, with $d$, $g$ and $u_0$ given and $\partial \Omega = \Gamma_D \cup \Gamma_N$, $\Gamma_D \cap \Gamma_N = \emptyset$. 
The \emph{Reynolds} number $Re$ depends on the kinematic viscosity $\nu$ through the characteristic velocity $U$ and characteristic length $L$ as 

\begin{equation}
 Re = \frac{UL}{\nu}.   
 \label{Hess:Re_def}
\end{equation}

As an example of the inner workings of Nektar++, the computation of steady states is explained in detail, which will used for the model reduction in the next section.
The steady states are solutions where $\frac{\partial u}{\partial t} = 0$ holds.

For the ROM in the next section, a parametric variation of the viscosity $\nu$ will be assumed.
This corresponds to a variable Reynolds number via the relation \eqref{Hess:Re_def} and allows to compute the flow for various Reynolds numbers.
To compute steady states for varying viscosity $\nu$, a solution $u(\nu_0)$ for a parameter value $\nu_0$ is used as an initial guess within a fixed point iteration to obtain the steady state solution $u(\nu_1)$ at a  parameter value $\nu_1$. 
This is repeated for $\nu_2$ and so on. In this way, the parameter range of interest can be explored iteratively with a fixed point iteration.

The \emph{Oseen}-iteration is a secant modulus fixed-point iteration with a linear rate of convergence. Given the current iterate (or initial condition) $u^k$, the linear system
\begin{eqnarray}
 -\nu \Delta u + (u^k \cdot \nabla) u + \nabla p &=& f  \text{ in } \Omega, \label{Hess:eq_Oseen_main}\\
\nabla \cdot u &=& 0   \text{ in } \Omega, \\
u &=& d \quad \text{ on } \Gamma_D, \\
\nabla u \cdot n &=& g \quad \text{ on } \Gamma_N, 
\end{eqnarray}
\noindent is solved for the next iterate $u^{k+1} = u$. 
A usual stopping criterion is that the relative change between iterates in the $L^2$ or $H^1$ norm falls below a given tolerance.
The initial solution $u^0(\nu_0)$ is computed by time-advancement of \eqref{Hess:NSE0}--\eqref{Hess:NSE1} from zero initial conditions at a parameter value $\nu_0$. From this starting point, solutions on the whole parameter domain can be found.

The discretized system solved in each step of the \emph{Oseen}-iteration is decomposed as \eqref{Hess:fully_expanded} as
\begin{eqnarray}
\left[
\begin{array}{ccc}
A & -D^T_{bnd} & B  \\
-D_{bnd} & 0 & -D_{int} \\
\tilde{B}^T & -D^T_{int} & C
\end{array}
\right]
\left[
\begin{array}{ccc}
v_{bnd} \\
p \\
v_{int} 
\end{array}
\right]
&=
\left[
\begin{array}{ccc}
f_{bnd} \\
0 \\
f_{int}
\end{array}
\right]
\label{Hess:fully_expanded}
\end{eqnarray}
\noindent where $v_{bnd}$ and $v_{int}$ denote velocity degrees of freedom on the boundary and in the interior, respectively.
The forcing terms $f_{bnd}$ and $f_{int}$ refer to the boundary and interior, respectively.
The matrix $A$ assembles the boundary-boundary terms, $B$ the boundary-interior terms, $\tilde{B}$ the interior-boundary terms and
$C$ assembles the interior-interior terms of elemental velocity \emph{ansatz} functions. 
In a \emph{Stokes} system, it holds that $B = \tilde{B}^T$, but this is not the case in the \emph{Oseen} equation, since 
the linearization term $(u^k \cdot \nabla) u$ is present in \eqref{Hess:eq_Oseen_main}.
The matrices $D_{bnd}$ and $D_{int}$ provide the pressure-velocity boundary and pressure-velocity interior couplings, respectively.

The linear system \eqref{Hess:fully_expanded} is assembled in local degrees of freedom, leading to block matrices $A, B, \tilde{B}, C, D_{bnd}$ and $D_{int}$, with
each block referring to one spectral element. Thus, the system is singular in this form.
To solve the system, the local degrees of freedom can be gathered into the global degrees of freedom, but here a multi-level static condensation is employed. See also the documentation and source code of \emph{Nektar++}.

\subsection{Reduced order modeling for continuous Galerkin methods}
The reduced order model (ROM) needs to approximate the full order solutions accurately over the parameter domain of interest.
The reduced basis (RB) model reduction uses a projection onto a low order space of snapshot solutions (i.e., full order solutions) and 
an offline-online decomposition to facilitate computational efficiency.
A set of snapshots is computed 
over a coarse sample of the parameter domain 
and used to define a projection space $U$ of size $N$ using the standard solver provided by Nektar++. 
The POD computes a singular value decomposition (SVD) of the snapshot solutions to $99.9\%$ of the most dominant POD modes.
This defines a projection matrix $U \in \mathbb{R}^{N_\delta \times N}$.
The software package \textit{Eigen} is used in ITHACA-SEM to compute the SVD.

The offline-online decomposition allows for fast input-output solves, because they are independent of the original model size $N_\delta$.
It is an important part of efficient reduced order modeling, but since the static condensation includes 
the inversion of the parameter-dependent matrix $C$, the projection is applied to the system  \eqref{Hess:fully_expanded}. Alternatively, also some degrees of freedom can be gathered, see \cite{HessRozza2019}.
In the offline phase, snapshot solutions have been gathered over the parameter domain, which now serve as a projection space to define the reduced order setting. 
To have fast reduced order solves, the offline-online decomposition expands \eqref{Hess:fully_expanded} in the parameter of interest and 
stores the parameter independent projections as small-sized matrices of the order $N \times N$.
Since during the \emph{Oseen}-iteration each matrix is dependent on the previous iterate, the submatrices corresponding to each basis function are assembled and then formed online with the reduced basis coordinate representation of the current iterate. 
This is analogous to reduced order assembly of the nonlinear term in the \emph{Navier-Stokes} equations.
For more details and applications, see \cite{10.1007/978-3-030-39647-3_45}.

Particular care must be taken regarding the different levels of accuracy of a function within Nektar++.
There exist three levels, the global degrees of freedom, the local degrees of freedom, where degrees of freedom on the boundaries between spectral elements are present multiple times and the physical degrees of freedom, which correspond to quadrature points. The nonlinear terms need to be evaluated on the level of physical degrees of freedom for the increased accuracy but the projection takes place in the local degrees of freedom. This correspondingly requires two separate sets of projection bases.

\subsubsection{A numerical example}

\begin{figure}
 \includegraphics[scale=.2]{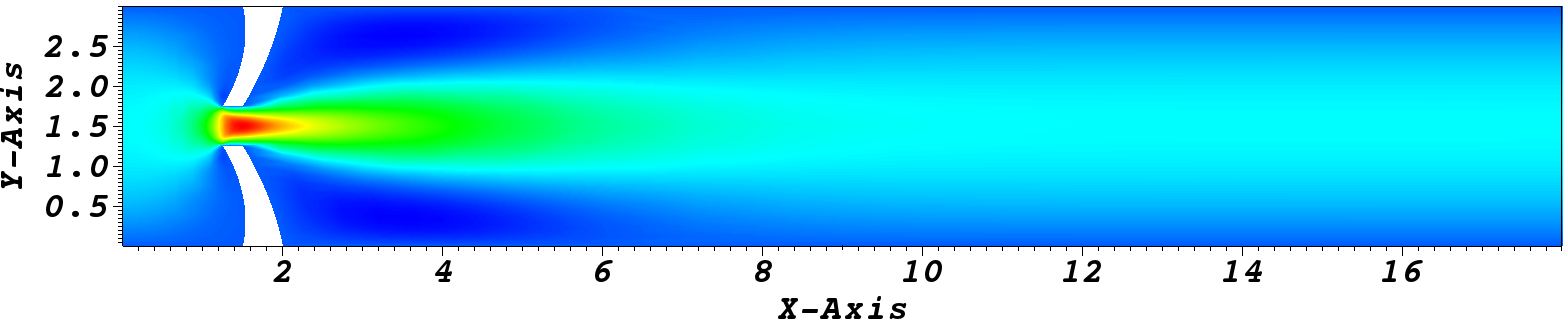} $\quad$
 \includegraphics[scale=.3]{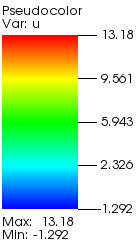} \\
 \includegraphics[scale=.2]{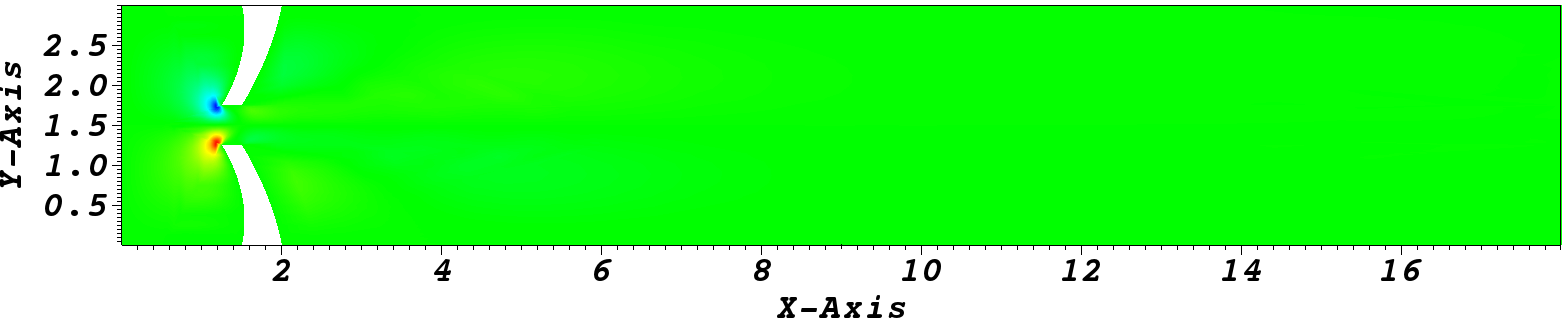} $\quad$
 \includegraphics[scale=.3]{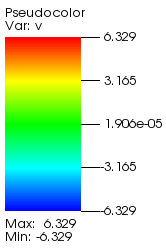} 
 \caption{Full order, steady-state solution for $\nu = 0.15$: velocity in x-direction (top) and y-direction (bottom).}
 \label{Hess:nu10}
\end{figure}

\begin{figure}
 \includegraphics[scale=.2]{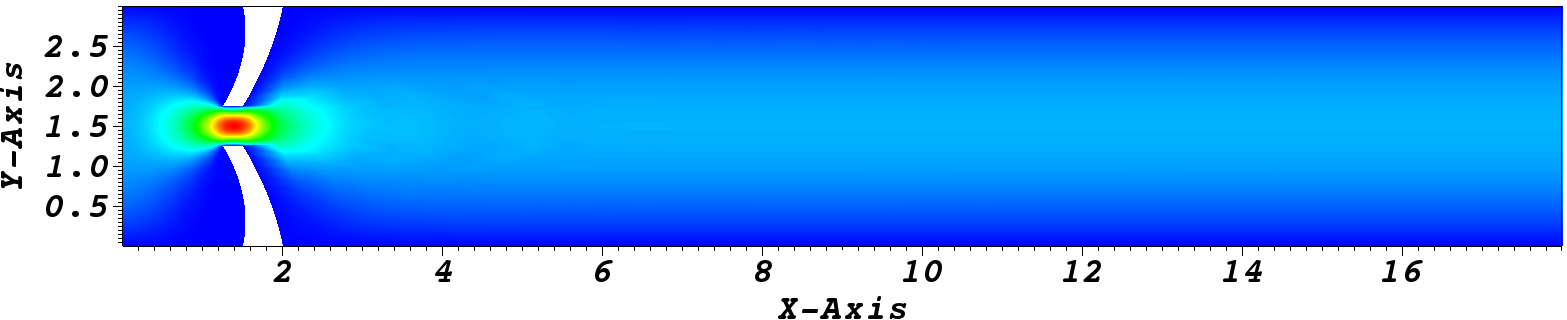} $\quad$
 \includegraphics[scale=.3]{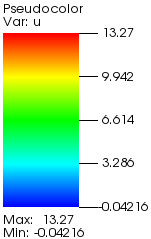} \\
 \includegraphics[scale=.2]{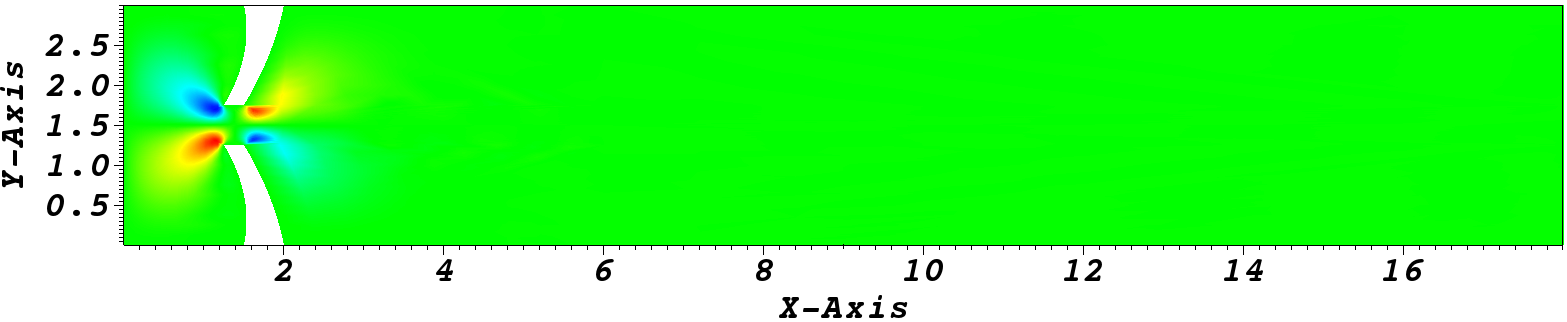} $\quad$
 \includegraphics[scale=.3]{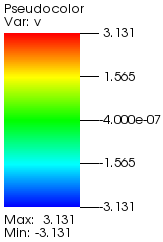} 
 \caption{Full order, steady-state solution for $\nu = 10$: velocity in x-direction (top) and y-direction (bottom).}
 \label{Hess:nu0p15}
\end{figure}

Consider a variable kinematic viscosity in the interval $\nu \in [ 0.15 ; 10 ]$ for a channel flow as depicted in Fig.~\ref{Hess:nu10} and Fig.~\ref{Hess:nu0p15}. A parabolic inflow profile is prescribed on the left wall at $y=0$, a natural outflow boundary is prescribed at $y=8$ and the other walls are no-slip boundaries. 
The polynomial order of the ansatz functions is chosen as 11. The relative  error between full order solutions upon increasing the polynomial order is about $0.05\%$.
The Fig.~\ref{Hess:nu10} and Fig.~\ref{Hess:nu0p15} show the extreme cases considered here for a very small viscosity and very large viscosity.

There is a bifurcation occurring nearby for some $\nu < \nu^{\text{min}} = 0.15$ as investigated in \cite{PITTON2017534}. Since the bifurcation point represents a singularity, the convergence speed of the fixed-point iteration is very low close to the bifurcation. The tolerance for the fixed-point iteration has thus been set $1\mathrm{e}{-4}$ for a change among two iterates.

\begin{figure}
 \includegraphics[scale=.8]{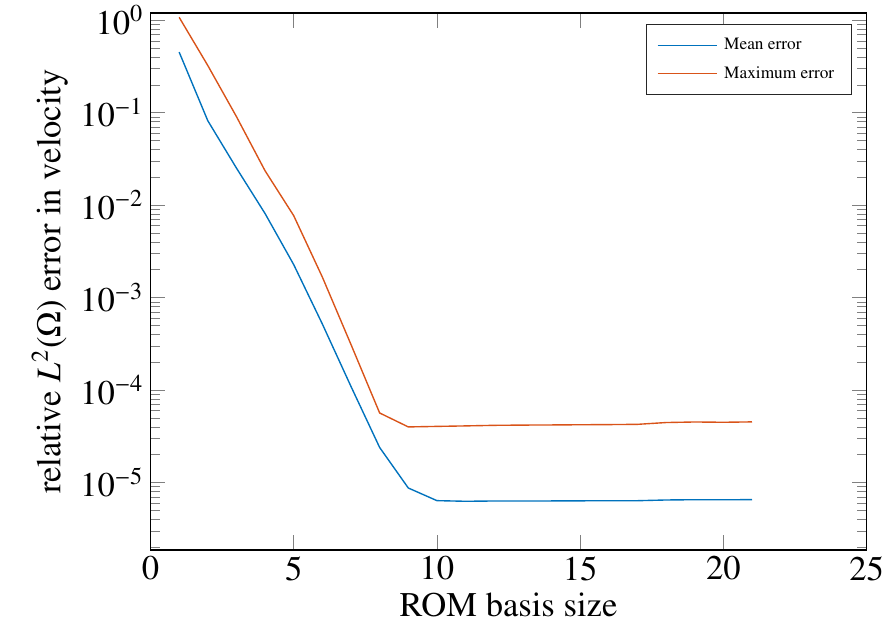}
 \caption{Mean and maximum relative $L^2(\Omega)$ error in the velocity with increasing basis size.}
 \label{Hess:l2_error_decay}
\end{figure}

Sample solutions are computed at $22$ values of $\nu$ in the interval of interest $[ 0.15 ; 10 ]$. The POD is computed and the 
POD energy reaches a threshold of $99\%$ with 2 modes and a threshold of $99.99\%$ with 6 modes.
Ideally, the exponential decay in POD energy should translate into an exponential decay in relative approximation error.
Fig.~\ref{Hess:l2_error_decay} shows the mean and maximum relative $L^2(\Omega)$ error in the velocity as an increasing number of POD modes is considered. A mean relative error of $1\%$ is reached with 5 basis functions and a maximum relative error of $1\%$ is reached with 6 basis functions. 

The error plateaus at about $1\mathrm{e}{-4}$, which is because the fixed-point iteration is not more accurate than that. When the problem is such that the fixed-point iteration can work with a higher accuracy, then this will also be recovered in the ROM.
The offline time took 46s and the online evaluation at all parameter values and over all basis sizes took 1s on a workstation with an i7-6700.

\subsubsection{Parametric variation in geometry}

Numerical examples where a parametrized geometry is considered can be found in \cite{HessQuainiRozza2020} and \cite{10.1007/978-3-030-39647-3_45}.
In case of an affine parameter dependency, this parameter dependency can be made explicit 
in the system matrix as 

\begin{equation}
\label{Hess:affine_full_dim}
 A(\mu) x = \left( \sum_{i=1}^Q  \Theta_i(\mu) A_i \right) x =  b,
\end{equation}

\noindent given a parameter vector $\mu$ with scalar functions $\Theta_i$ and parameter independent matrices $A_i$. The matrix $A$ corresponds to the matrix  \eqref{Hess:fully_expanded} in the full order assembly.
The precise transformations to arrive at the form \eqref{Hess:affine_full_dim} can be found in \cite{QuarteroniRozza2007affineGEO}. ITHACA-SEM allows the user to specify the form
 \eqref{Hess:affine_full_dim} in a header file which is then used to achieve the offline-online decomposition by exploiting
 
 \begin{equation}
\label{Hess:affine_reduced_dim}
 V^T A(\mu) V x_r = \left( \sum_{i=1}^Q  \Theta_i(\mu) V^T A_i V \right)  x_r =  V^T b ,
\end{equation}

\noindent given a projection matrix $V$. The terms $V^T A_i V$ are computed in the offline phase
and are then available online as small scale matrices, which form the reduced order system by assembling the sum for a parameter of interest $\mu$. The approximation to the full-order solution $x$ is then given by 

\begin{equation}
    x \approx V x_r.
\label{Hess:approx_x}
\end{equation}

In case of a parametric variation in geometry which is non-affine, the \emph{empirical interpolation method} (EIM) is usually used, see \cite{eim1}. The EIM approximates an affine form, such as \eqref{Hess:affine_full_dim}. To do that, a few degrees of freedom are identified, which allow to recover the whole system matrix in an affine form. Its speed-up depends on the ability to quickly compute a few entries of the system matrix. This is typically the case if the \emph{ansatz} functions have a small support such as in a finite element method (FEM). But in the SEM the \emph{ansatz} functions have a much larger support, such that this speed-up is not possible and the applicability of the EIM is limited. This was investigated in \cite{HessQuainiRozza2020}.

\section{Discontinuous Galerkin: compressible flow simulations}
\label{sec:dg}

\subsection{Overview}
The DG-based solver implemented in Nektar++ is used for the prediction of compressible flows in which density variations are not negligible.
This implies that the zero-divergence condition of the velocity fields no longer holds; the resulting set of equations follows:
\begin{eqnarray}
\frac{\partial \rho}{\partial t} +  \nabla \cdot (\rho u) &=& 0 \\ \nonumber 
\frac{\partial u}{\partial t} +   \nabla(\rho u \times u) &=& - \nabla p + \nabla\cdot \sigma+ \rho f, \\
\frac{\partial E}{\partial t} + \nabla\cdot(u E) + \nabla\cdot(u p) &=&  \nabla \cdot(u\cdot \sigma) + k\nabla T \quad, \nonumber
\label{DG:NSEcmpr}
\end{eqnarray}
where $\rho$ is the density, $E$ is the total energy, and $T$ the temperature. Moreover $\sigma$ indicates the viscous dissipation tensor:
\begin{equation}
\sigma = 2\nu\rho \left( S_{ij}+\frac{1}{3}  \nabla \cdot u \delta_{ij }\right) \quad,
\end{equation}
being the strain tensor $S_{ij}$ defined as the symmetric part of the gradient of the velocity $u$. To close the problem, proper initial and boundary conditions must be provided:
\begin{eqnarray}
U_i &=& d \quad \text{ on } \Gamma_D \times (0, T), \\ \nonumber
\nabla U_i \cdot n &=& g \quad \text{ on } \Gamma_N \times (0, T), \\
U_i &=& U_{i0} \quad \text{ in } \Omega \times 0,  \nonumber
\label{DG:NSE_boundaryCond} 
\end{eqnarray}
being $U={U_1, U_2, U_3} = {\rho, u, T}$.

In detail the DG scheme is a high order method that, differently from the continuous Galerkin approach, does not impose the continuity of the solution between contiguous elements.
In order to keep the consistency of the solution, proper flux must be exchanged and in general flux splitting schemes are devised by the ones developed in the Finite Volume Method framework.
In detail, in this work the HLLC flux splitting scheme was used. The HLLC differs from the more common HLL scheme because it is able to reconstruct the contact surface of a Riemann problem.
For further details, one can refer to~\cite{toro1994restoration}.
Time integration was performed with an explicit fourth order Runge-Kutta marching scheme, chosen because of its low dissipative properties.
A more comprehensive description of the full order method implemented in Nektar++ can be found in~\cite{cantwell2015nektar++}.

\subsection{Reduced order modeling for discontinuous Galerkin methods}
As in the previous case, the ROM for compressible flows is structured according to an online/offline paradigm.
During the offline phase, once that $N$ snapshots have been generated with the aid of the full order solver implemented in Nektar++ and gathered in a database, the dominant dynamics are extracted through a SVD procedure.
The computational speed up of the ROMs relies on the fact that only the first  $N_m < N$ most energetic modes are retained when the building of the bases is carried on, being $N$ the number of degree of freedom of the full order solver.
The equations \eqref{DG:NSEcmpr} are linearized and the operators are projected on the reduced basis once that a suitable inner product has been identified.
ROMs for compressible flows are challenging since an energy-based inner product is not defined and the evolution of ROMs does not satisfy the energy equation.
In this work a stabilization matrix inspired by the symmetry matrix introduced by Barone et al.~\cite{barone2009stable} is used for recovering part of the stability of the reduced order system.

The introduction of the stabilization matrix aims to improve, beyond the stability, also the controllability of the system, e.g.\ moving the unstable (positive) eigenvectors closer to zero. This has important consequences when Eigenvalues Replacement (ER) techniques are employed for stabilizing the system since smaller spaces have to be explored. 
The first method based on Eigenvalues Replacement can be found in Balajewicz et al.~\cite{balajewicz2016minimal}, who propose to achieve the stability of the reduced system through a minimal subspace rotation, thus by introducing viscosity.
In this work the ER is based on a Swarm Particle Optimization in order to find the combination of eigenvalues which minimize the distance between the predicted coefficients $a_{ROM}$ and the ones $a_{FOM}$ obtained by projecting the snapshots onto the reduced basis \cite{rezaian2018eigenvalue}.
Moreover one requires that the total power of the system W decreases with time. Hence:
\begin{equation}
min_{\lambda} \sum_i^{N_s} || a_{ROM}^k - a_{FOM}^k ||_{2}^{2}  \quad s.t. \quad W(t) < 0 \quad,
\label{DG:opt}
\end{equation}
being $N_s$ the number of snapshots and the total power computed as $E \approx a_{ROM}^2 $.
To include the constraints related to the total power, Equation~\ref{DG:opt} can be rewritten as follows:
\begin{equation}
min_{\lambda} \sum_i^{N_s} || a_{ROM}^k - a_{FOM}^k ||_{2}^{2}  +c1( \alpha + c_2) \quad,
\label{DG:opt2}
\end{equation}
where $c_1$ and $c2$ are constants, $\alpha$ is the coefficient of the linear regression which gives the best fit of $W(t) = \alpha  t + \beta$.
In general $c_1$ is a constant which represent a penalization term and it has to be determined each time, while $c_2$ is a small positive number, that in the present work is evaluated as:
\begin{equation}
c_2 = max(10^{-5}, \alpha_{FOM}) \quad,
\end{equation}
where $\alpha_{FOM}$ is the angular coefficient of the linear regression which best approximates the total power computed starting from the coefficients obtained from the full order snapshots.
The time integration is performed with an explicit fourth order Runge Kutta scheme in order to be consistent with the full order solver.

\subsection{Results}
\begin{figure}
    \centering
    \includegraphics[width=.7\textwidth]{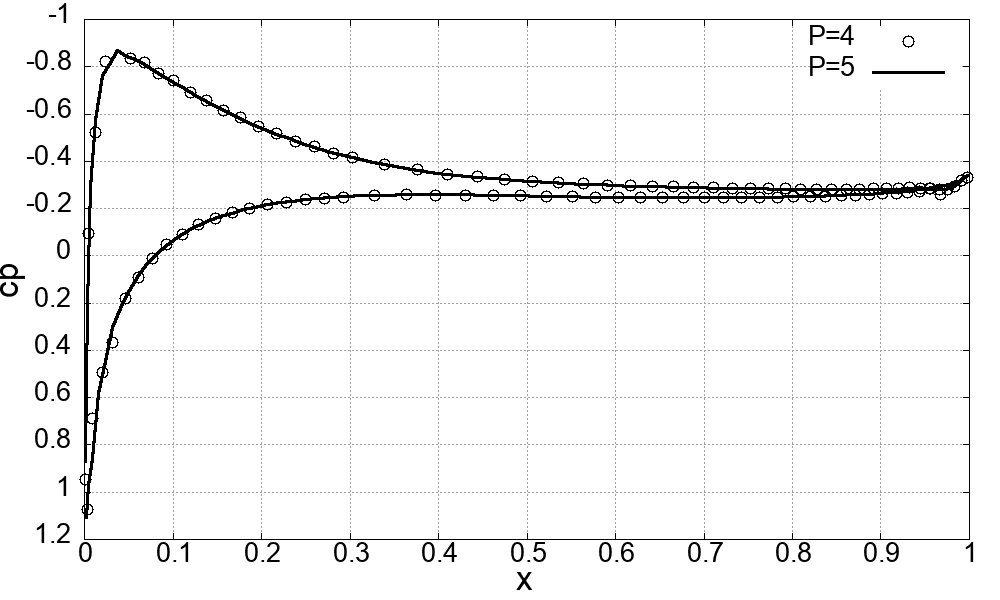}
    \caption{Time-averaged pressure coefficients over the airfoil obtained with two different polynomial orders: $P=4$ (empty dots), $P=5$ (solid line).}
    \label{fig:DG_4vs5}
\end{figure}

\begin{figure}
\centering
\begin{subfigure}[b]{0.495\textwidth}
 \includegraphics[width=1\textwidth]{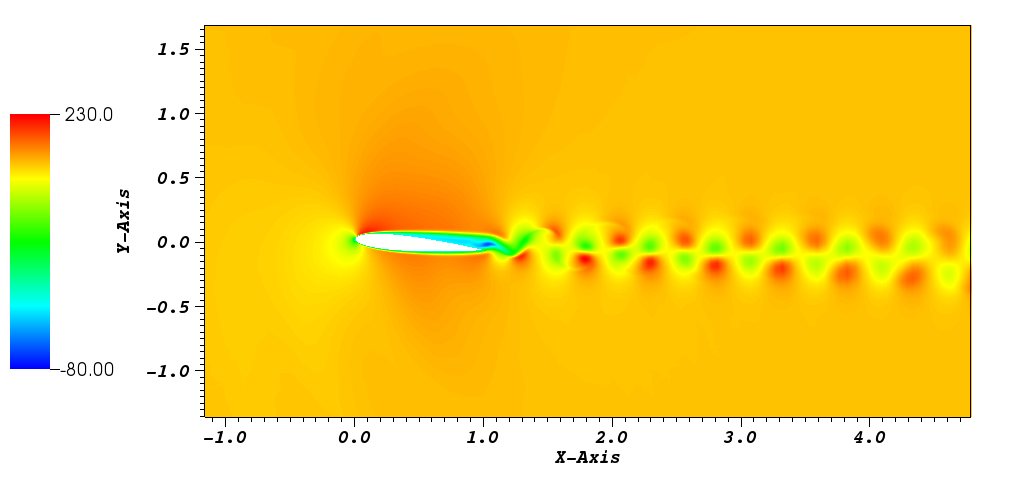}
\end{subfigure}
     \hfill
\begin{subfigure}[b]{0.495\textwidth}
 \includegraphics[width=1\textwidth]{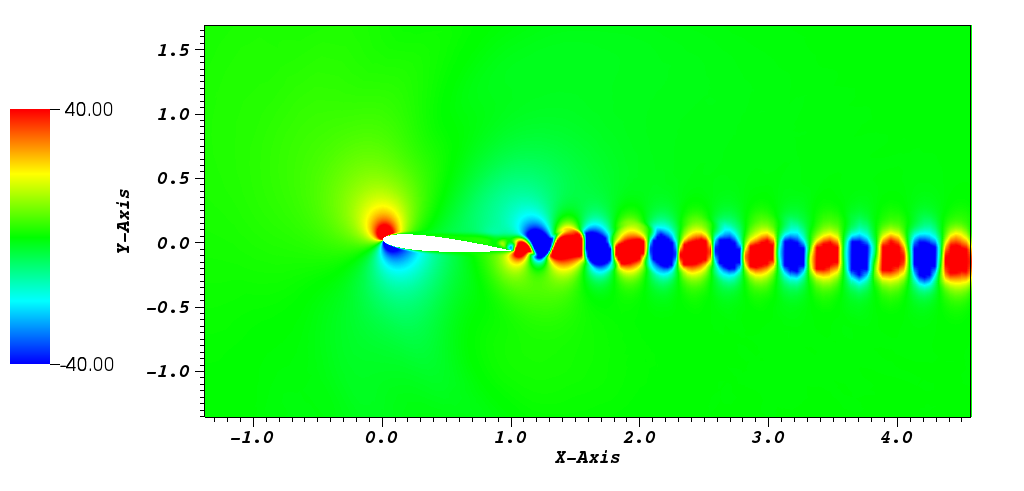} 
\end{subfigure}
\caption{Flow field for Re=7500: velocity along the x-axis (left) and y-axis (right).}
\label{DG:refCond}
\end{figure}

A NACA 0012 airfoil exposed to a flow field with an angle of attack of 5 degrees is considered as a test case for applying the reduction strategy described in the previous sections.
The unperturbed flow field is characterized by a Mach number of 0.5 and standard sea level air conditions are used for initializing the flow variables.
The full order snapshots have been obtained using a two-dimensional computational grid whose main dimensions were 40 times the airfoil chord length and it was composed by approximately 4000 discrete elements. Moreover, a polynomial order $P$ equal to four was used for the simulation. Figure \ref{fig:DG_4vs5} shows the distributions of pressure coefficients along the x-axis obtained by using two distinct polynomial orders (four and five); the two solutions do not differ significantly, thus indicating the convergence of the numerical model for $P=4$.
The resulting distribution of velocity components along the x-axis and y-axis obtained with the full order solver for a Reynolds number equal to 7500 is shown in Figure~\ref{DG:refCond}.

\begin{figure}
\centering
\begin{subfigure}[b]{0.495\textwidth}
 \includegraphics[width=1\textwidth]{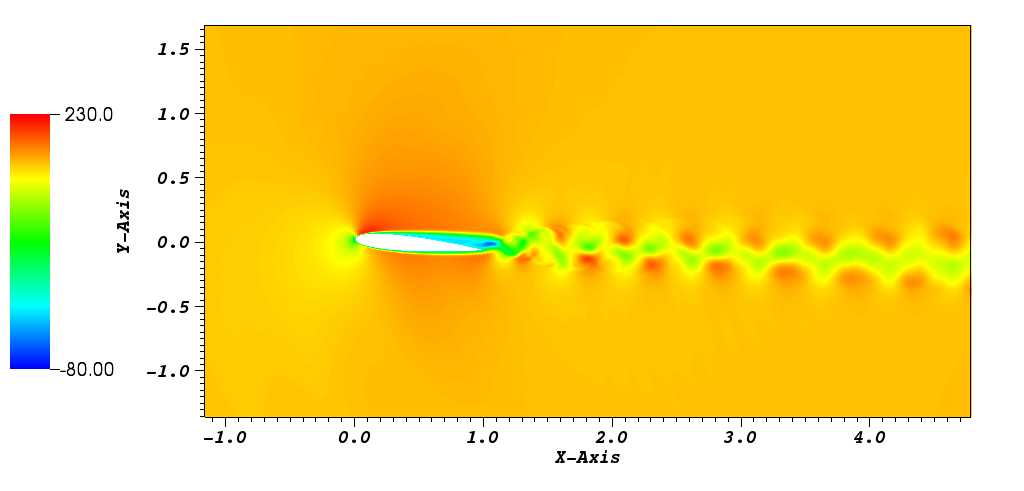}
\end{subfigure}
     \hfill
\begin{subfigure}[b]{0.495\textwidth}
 \includegraphics[width=1\textwidth]{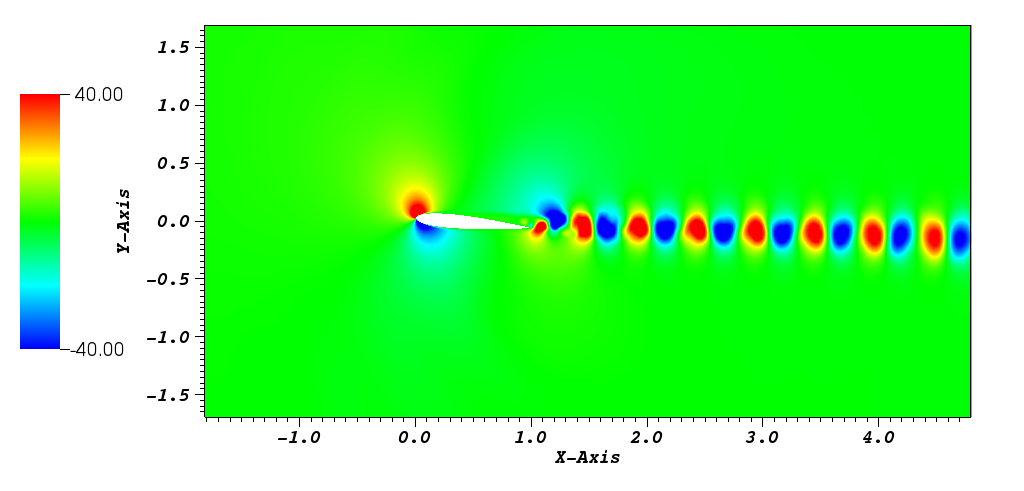} 
\end{subfigure}
\begin{subfigure}[b]{0.49\textwidth}
 \includegraphics[width=1\textwidth]{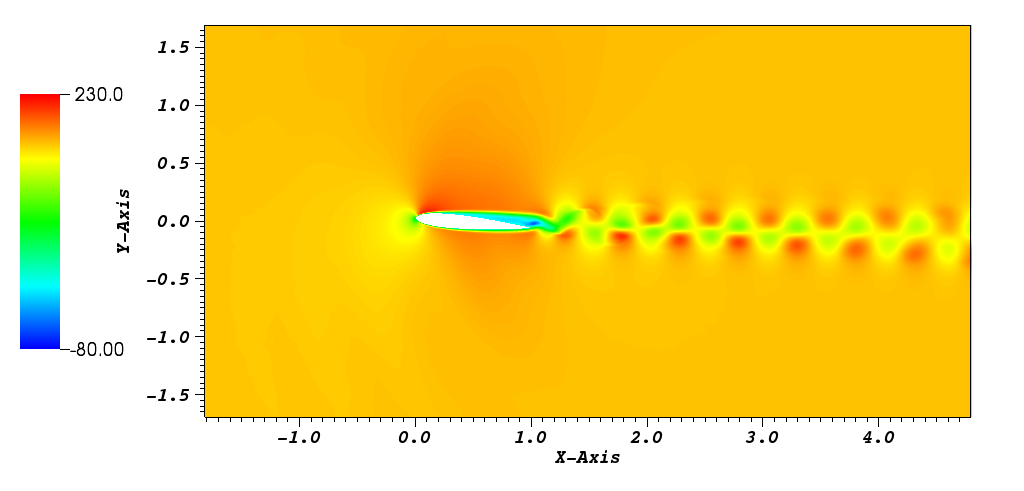}
\end{subfigure}
     \hfill
\begin{subfigure}[b]{0.49\textwidth}
 \includegraphics[width=1\textwidth]{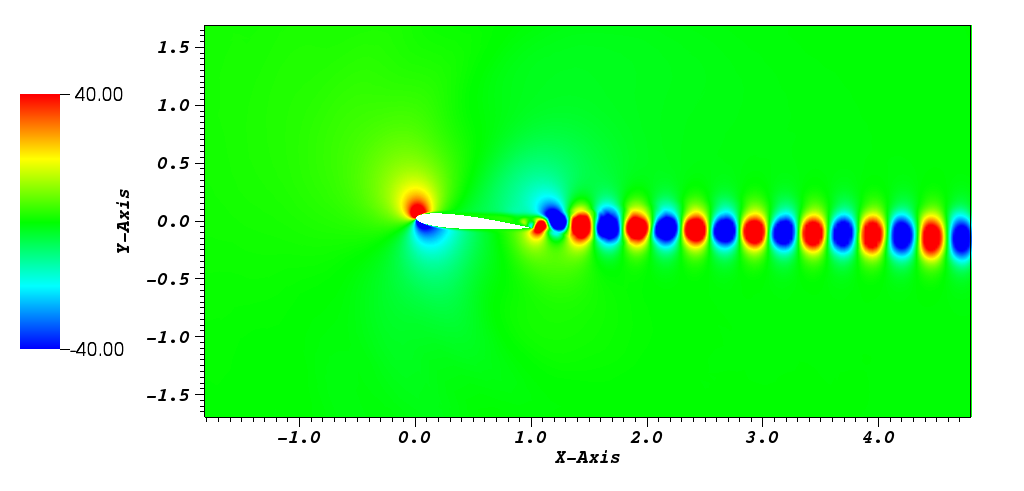} 
\end{subfigure}
\caption{Velocity along the x-axis (left) and y-axis (right): comparison between a simple stabilization strategy (top) and one reinforced with an eigenvalue replacement (bottom)}
\label{DG:reduction}
\end{figure}

For each flow realization, 60 snapshots were saved and 8 modes were generated.
Two distinct reconstruction strategies were adopted: in the first case the linearized governing equations were projected on the reduced basis and the reduced matrices were computed; in the second case the eigenvalues of these matrices were further optimized with the aid of the PSO-strategy previously described.
Reconstructed flow fields for Re=7500 are reported in Figure~\ref{DG:reduction}; as one can observe the additional step contribute to enhance the stability and this can be seen in particular for the u component of the velocity vector.
 The computational effort required by the online phase is approximately 28 seconds, on the same machine the FOM solution requires approximately 5 hours, given the small time step needed to guarantee the stability of the solution.

\section{Conclusions}
\label{sec:conclusions}
In this paper, we presented the first steps towards intrusive ROM capabilities for SEM. In particular, we highlighted the implementation of offline-online ROM decomposition strategies for both CG and DG, in the context of incompressible and compressible flow problems, respectively. These steps constitutes the basis for further ROM developments for SEM, that is of particular interest for high-fidelity flow simulations in complex geometries, such as explored for finite volume methods in \cite{Stabile20202655, HijaziStabileMolaRozza2020}. In detail next steps will be focused on the development of non-intrusive Reduced Order Models based on Machine Learning techniques and on the stabilization on non-linear compressible ROMs. Indeed, non-intrusive methods present some advantages with respect to intrusive ones, among the others they are naturally stable and more user friendly to deal with.


%
\begin{acknowledgement}
We acknowledge the support by European
Union Funding for Research and Innovation - Horizon 2020 Program - in the framework of European Research
Council Executive Agency: Consolidator Grant H2020 ERC CoG 2015 AROMA-CFD project 681447 “Advanced
Reduced Order Methods with Applications in Computational Fluid Dynamics” (P.I.Prof. Gianluigi Rozza).
This work was partially supported by MIUR (Italian ministry for university and research) through FARE-X-AROMA-CFD project, P.I. Prof. Gianluigi Rozza. Gianmarco Mengaldo acknowledges support from NUS startup grant R-265-000-A36-133.
\end{acknowledgement}

\bibliographystyle{spmpsci}
\bibliography{references}

\begin{thebibliography}{10}
\providecommand{\url}[1]{{#1}}
\providecommand{\urlprefix}{URL }
\expandafter\ifx\csname urlstyle\endcsname\relax
  \providecommand{\doi}[1]{DOI~\discretionary{}{}{}#1}\else
  \providecommand{\doi}{DOI~\discretionary{}{}{}\begingroup
  \urlstyle{rm}\Url}\fi

\bibitem{balajewicz2016minimal}
Balajewicz, M., Tezaur, I., Dowell, E.: Minimal subspace rotation on the
  {S}tiefel manifold for stabilization and enhancement of projection-based
  reduced order models for the compressible {N}avier--{S}tokes equations.
\newblock Journal of Computational Physics \textbf{321}, 224--241 (2016)

\bibitem{barone2009stable}
Barone, M.F., Kalashnikova, I., Segalman, D.J., Thornquist, H.K.: Stable
  {G}alerkin reduced order models for linearized compressible flow.
\newblock Journal of Computational Physics \textbf{228}(6), 1932--1946 (2009)

\bibitem{cantwell2015nektar++}
Cantwell, C.D., Moxey, D., Comerford, A., Bolis, A., Rocco, G., Mengaldo, G.,
  De~Grazia, D., Yakovlev, S., Lombard, J.E., Ekelschot, D., et~al.: Nektar++:
  An open-source spectral/hp element framework.
\newblock Computer physics communications \textbf{192}, 205--219 (2015)

\bibitem{deville_fischer_mund_2002}
Deville, M.O., Fischer, P.F., Mund, E.H.: High-Order Methods for Incompressible
  Fluid Flow.
\newblock Cambridge Monographs on Applied and Computational Mathematics.
  Cambridge University Press (2002).
\newblock \doi{10.1017/CBO9780511546792}

\bibitem{fernandez2019non}
Fernandez, P., Moura, R.C., Mengaldo, G., Peraire, J.: Non-modal analysis of
  spectral element methods: Towards accurate and robust large-eddy simulations.
\newblock Computer Methods in Applied Mechanics and Engineering \textbf{346},
  43--62 (2019)

\bibitem{eim1}
Grepl, M.A., Maday, Y., Nguyen, N.C., Patera, A.T.: Efficient reduced-basis
  treatment of nonaffine and nonlinear partial differential equations.
\newblock ESAIM: Mathematical Modelling and Numerical Analysis \textbf{41}(3),
  575--605 (2007).
\newblock \doi{10.1051/m2an:2007031}

\bibitem{HessQuainiRozza2020}
Hess, M.W., Quaini, A., Rozza, G.: Reduced basis model order reduction for
  {N}avier-{S}tokes equations in domains with walls of varying curvature.
\newblock International Journal of Computational Fluid Dynamics \textbf{34}(2),
  119--126 (2020).
\newblock \doi{10.1080/10618562.2019.1645328}

\bibitem{10.1007/978-3-030-39647-3_45}
Hess, M.W., Quaini, A., Rozza, G.: A spectral element reduced basis method for
  {N}avier--{S}tokes equations with geometric variations.
\newblock In: S.J. Sherwin, D.~Moxey, J.~Peir{\'o}, P.E. Vincent, C.~Schwab
  (eds.) Spectral and High Order Methods for Partial Differential Equations
  ICOSAHOM 2018, pp. 561--571. Springer International Publishing, Cham (2020)

\bibitem{HessRozza2019}
Hess, M.W., Rozza, G.: A Spectral Element Reduced Basis Method in Parametric
  CFD, vol. 126, chap. A Spectral Element Reduced Basis Method in Parametric
  CFD, pp. 693--701.
\newblock Springer International Publishing (2019).
\newblock \doi{10.1007/978-3-319-96415-7_64}

\bibitem{hesthaven2015certified}
Hesthaven, J., Rozza, G., Stamm, B.: Certified Reduced Basis Methods for
  Parametrized Partial Differential Equations.
\newblock Springer Briefs in Mathematics (2015)

\bibitem{HijaziStabileMolaRozza2020}
Hijazi, S., Stabile, G., Mola, A., Rozza, G.: Data-driven {POD}-{G}alerkin
  reduced order model for turbulent flows.
\newblock Journal of Computational Physics \textbf{416}, 109513 (2020).
\newblock \doi{10.1016/j.jcp.2020.109513}

\bibitem{karniadakis1999spectral}
Karniadakis, G., Sherwin, S.: Spectral/hp Element Methods for CFD.
\newblock Numerical mathematics and scientific computation. Oxford University
  Press (2005)

\bibitem{mengaldo2015discontinuous}
Mengaldo, G.: Discontinuous spectral/hp element methods: development, analysis
  and applications to compressible flows.
\newblock Imperial College London (2015)

\bibitem{mengaldo2018spatial}
Mengaldo, G., Moura, R., Giralda, B., Peir{\'o}, J., Sherwin, S.: Spatial
  eigensolution analysis of discontinuous galerkin schemes with practical
  insights for under-resolved computations and implicit les.
\newblock Computers \& Fluids \textbf{169}, 349--364 (2018)

\bibitem{mengaldo2021industry}
Mengaldo, G., Moxey, D., Turner, M., Moura, R.C., Jassim, A., Taylor, M.,
  Peiro, J., Sherwin, S.: Industry-relevant implicit large-eddy simulation of a
  high-performance road car via spectral/hp element methods.
\newblock SIAM Review \textbf{63}(4), 723--755 (2021)

\bibitem{moura2020spatial}
Moura, R.C., Aman, M., Peir{\'o}, J., Sherwin, S.J.: Spatial eigenanalysis of
  spectral/hp continuous {G}alerkin schemes and their stabilisation via
  dg-mimicking spectral vanishing viscosity for high reynolds number flows.
\newblock Journal of Computational Physics \textbf{406}, 109112 (2020)

\bibitem{moura2020viscous}
Moura, R.C., Fernandez, P., Mengaldo, G., Sherwin, S.J.: Viscous diffusion
  effects in the eigenanalysis of (hybridisable) dg methods.
\newblock In: Spectral and High Order Methods for Partial Differential
  Equations ICOSAHOM 2018, pp. 371--382. Springer, Cham (2020)

\bibitem{moura2017setting}
Moura, R.C., Mengaldo, G., Peir{\'o}, J., Sherwin, S.J.: An les setting for
  dg-based implicit les with insights on dissipation and robustness.
\newblock In: Spectral and High Order Methods for Partial Differential
  Equations ICOSAHOM 2016, pp. 161--173. Springer, Cham (2017)

\bibitem{moura2017eddy}
Moura, R.C., Mengaldo, G., Peir{\'o}, J., Sherwin, S.J.: On the eddy-resolving
  capability of high-order discontinuous {G}alerkin approaches to implicit
  les/under-resolved dns of {E}uler turbulence.
\newblock Journal of Computational Physics \textbf{330}, 615--623 (2017)

\bibitem{moura2016eigensolution}
Moura, R.C., Sherwin, S.J., Peir{\'o}, J.: Eigensolution analysis of
  spectral/hp continuous galerkin approximations to advection--diffusion
  problems: Insights into spectral vanishing viscosity.
\newblock Journal of Computational Physics \textbf{307}, 401--422 (2016)

\bibitem{moxey2017towards}
Moxey, D., Cantwell, C., Mengaldo, G., Serson, D., Ekelschot, D., Peir{\'o},
  J., Sherwin, S., Kirby, R.: Towards p-adaptive spectral/hp element methods
  for modelling industrial flows.
\newblock In: Spectral and high order methods for partial differential
  equations icosahom 2016, pp. 63--79. Springer, Cham (2017)

\bibitem{moxey2020nektar++}
Moxey, D., Cantwell, C.D., Bao, Y., Cassinelli, A., Castiglioni, G., Chun, S.,
  Juda, E., Kazemi, E., Lackhove, K., Marcon, J., et~al.: Nektar++: Enhancing
  the capability and application of high-fidelity spectral/hp element methods.
\newblock Computer Physics Communications \textbf{249}, 107110 (2020)

\bibitem{PITTON2017534}
Pitton, G., Quaini, A., Rozza, G.: Computational reduction strategies for the
  detection of steady bifurcations in incompressible fluid-dynamics:
  Applications to {C}oanda effect in cardiology.
\newblock Journal of Computational Physics \textbf{344}, 534--557 (2017).
\newblock \doi{https://doi.org/10.1016/j.jcp.2017.05.010}.
\newblock
  \urlprefix\url{https://www.sciencedirect.com/science/article/pii/S0021999117303790}

\bibitem{QuarteroniRozza2007affineGEO}
Quarteroni, A., Rozza, G.: Numerical solution of parametrized
  {N}avier–{S}tokes equations by reduced basis methods.
\newblock Numerical Methods for Partial Differential Equations \textbf{23}(4),
  923--948 (2007).
\newblock \doi{https://doi.org/10.1002/num.20249}.
\newblock
  \urlprefix\url{https://onlinelibrary.wiley.com/doi/abs/10.1002/num.20249}

\bibitem{rezaian2018eigenvalue}
Rezaian, E., Wei, M.: Eigenvalue reassignment by particle swarm optimization
  toward stability and accuracy in nonlinear reduced-order models.
\newblock In: 2018 Fluid Dynamics Conference, p. 3095 (2018)

\bibitem{SchildersMOR}
Schilders, W.H., van~der Vorst, H.A., Rommes, J. (eds.): Model Order Reduction:
  Theory, Research Aspects and Applications.
\newblock Mathematics in Industry. Springer (2008).
\newblock \doi{https://doi.org/10.1007/978-3-540-78841-6}

\bibitem{Stabile20202655}
Stabile, G., Zancanaro, M., Rozza, G.: Efficient geometrical parametrization
  for finite-volume-based reduced order methods.
\newblock International Journal for Numerical Methods in Engineering
  \textbf{121}(12), 2655--2682 (2020).
\newblock \doi{10.1002/nme.6324}

\bibitem{tonicello2021fully}
Tonicello, N., Moura, R.C., Lodato, G., Mengaldo, G.: Fully-discrete spatial
  eigenanalysis of discontinuous spectral element methods: insights into
  well-resolved and under-resolved vortical flows.
\newblock arXiv preprint arXiv:2111.13891  (2021)

\bibitem{toro1994restoration}
Toro, E.F., Spruce, M., Speares, W.: Restoration of the contact surface in the
  hll-{R}iemann solver.
\newblock Shock waves \textbf{4}(1), 25--34 (1994)

\end{thebibliography}

\end{document}